%% file: Hilbert.tex
\input abbrev


\ref\BFH{J. de Boer, L. Feher, and A. Honecker, A class of $\cW$-algebras with infinitely generated classical limit, Nucl. Phys. B420 (1994), 409-445.}
\ref\Bor{R. Borcherds, Vertex operator algebras, Kac-Moody
algebras and the Monster, Proc. Nat. Acad. Sci. USA 83 (1986)
3068-3071.} 
\ref\BEHHH{R. Blumenhagen, W. Eholzer, A. Honecker, K. Hornfeck, and R. Hubel, Coset realizations of unifying $\cW$-algebras, Int. Jour. Mod. Phys. Lett. A10 (1995) 2367-2430.} 
\ref\DLM{C. Dong, H. Li, and G. Mason, Compact automorphism groups of vertex operator algebras, Internat. Math. Res. Notices 1996, no. 18, 913--921.}
\ref\EFH{W. Eholzer, L. Feher, and A. Honecker, Ghost systems: a vertex algebra point of view, Nuclear Phys. B 518 (1998), no. 3, 669--688.}
\ref\FKRW{E. Frenkel, V. Kac, A. Radul, and W. Wang, $\cW_{1+\infty}$ and $\cW(\gg\gl_n)$ with central charge $N$, Commun. Math. Phys. 170 (1995), 337-358.}
\ref\FHL{I.B. Frenkel, Y.Z. Huang, and J.
Lepowsky, On axiomatic approaches to vertex operator algebras and
modules, Mem. Amer. Math. Soc. 104 (1993), no. 494, viii+64.}
 \ref\FLM{ I.B. Frenkel, J. Lepowsky, and A. Meurman, Vertex
Operator Algebras and the Monster, Academic Press, New York, 1988.}
\ref\FMS{ D. Friedan, E. Martinec, and S. Shenker, Conformal invariance, supersymmetry and string theory, Nucl. Phys. B271 (1986) 93-165.} 
\ref\KP{V. Kac and D. Peterson, Spin and wedge representations of infinite-dimensional Lie algebras and groups Proc. Natl. Acad. Sci. USA 78 (1981), 3308-3312.}
\ref\KR{V. Kac and A. Radul, Representation theory of the vertex algebra $\cW_{1+\infty}$, Transf. Groups, Vol 1 (1996) 41-70.}
\ref\KWY{V. Kac, W. Wang, and C. Yan, Quasifinite representations of classical Lie subalgebras of $\cW_{1+\infty}$, Advances in Mathematics, vol. 139 (1), (1998) 59--140.}
 \ref\Kac{V. Kac, Vertex Algebras for Beginners, AMS Univ. Lecture Series Vol. 10, 2nd corrected ed., 2001.}
\ref\Li{H. Li, Local
systems of vertex operators, vertex superalgebras and modules,
J. Pure Appl. Algebra 109 (1996), no. 2, 143--195.}
\ref\Lii{ H. Li, Vertex algebras and vertex Poisson
algebras, Commun. Contemp. Math. 6 (2004) 61-110.}
 \ref\LZII{ B. Lian and G.J. Zuckerman,  Commutative quantum operator algebras, J. Pure Appl. Algebra 100 (1995) no. 1-3, 117-139.}
\ref\LL{B. Lian and A. Linshaw, Howe pairs in the theory of vertex algebras, J. Algebra 317, 111-152 (2007).}
\ref\LI{A. Linshaw, Invariant chiral differential operators and the $\cW_3$-algebra, J. Pure Appl. Algebra 213 (2009), 632-648.}
\ref\LII{A. Linshaw, Invariant theory and the $\cW_{1+\infty}$ algebra with negative integral central charge, arXiv:0811.4067.}
 \ref\MSV{F. Malikov, V. Schectman, and A. Vaintrob, Chiral de Rham complex, Comm. Math. Phys, 204, (1999) 439-473.}
\ref\WI{W. Wang, $\cW_{1+\infty}$ algebra, $\cW_3$ algebra, and Friedan-Martinec-Shenker bosonization, Commun. Math. Phys. 195 (1998), 95--111.}
\ref\WII{W. Wang, Classification of irreducible modules of $\cW_3$ with $c=-2$, Commun. Math. Phys. 195 (1998), 113--128.}
\ref\WIII{W. Wang, Duality in infinite-dimensional Fock representations, Commun. Contemp. Math. 1 (1999), 155--199.}
\ref\WIV{W. Wang, Dual pairs and infinite-dimensional Lie algebras, Contemp. Math. 248 (1999) 453--469.}
\ref\Zh{Y. Zhu, Modular invariants of characters of vertex operators, J. Amer. Math. Soc. 9 (1996) 237-302.}
\ref\We{H. Weyl, The Classical Groups: Their Invariants and Representations, Princeton University Press, 1946.}

\centerline{\titlefont A Hilbert theorem for vertex algebras}

\bs
\centerline{Andrew R. Linshaw}
\bs

\baselineskip=13pt plus 2pt minus 2pt

ABSTRACT. Given a simple vertex algebra $\cA$ and a reductive group $G$ of automorphisms of $\cA$, the invariant subalgebra $\cA^G$ is strongly finitely generated in most examples where its structure is known. This phenomenon is subtle, and is generally not true of the classical limit of $\cA^G$, which often requires infinitely many generators and infinitely many relations to describe. Using tools from classical invariant theory, together with recent results on the structure of the $\cW_{1+\infty}$ algebra, we establish the strong finite generation of a large family of invariant subalgebras of $\beta\gamma$-systems, $bc$-systems, and $bc\beta\gamma$-systems.

\baselineskip=15pt plus 2pt minus 1pt
\parskip=\baselineskip


\footnote{}{AMS Subject Classification: 17B69, 17B65, 13A50} 


\newsec{Introduction}
Let $G$ be a reductive group, and let $V$ be a finite-dimensional, linear representation of $G$. Throughout this paper, our base field will be $\C$. A basic problem in invariant theory is to describe the ring $\cO(V)^G$ of invariant polynomial functions on $V$. By a fundamental theorem of Hilbert, $\cO(V)^G$ is always finitely generated. The proof depends on the complete reducibility of $\cO(V)$ as a $G$-module, together with the fact that a summand of a finitely generated, graded commutative algebra, is itself finitely generated. It is also important to study $G$-invariant subalgebras of noncommutative rings $A$, such as universal enveloping algebras and Weyl algebras. Often, $A$ admits a filtration such that the associated graded algebra $gr(A)$ is commutative, and $G$ preserves the filtration. In this case, $(gr(A))^G\cong gr(A^G)$ as commutative rings, and information about $A^G$ such as a generating set can be obtained using commutative algebra. The ideal structure and the representation theory of $A^G$ may also be interesting, but unfortunately these cannot be reconstructed from $gr(A^G)$.

In this paper, we consider the invariant theory problem for {\it vertex algebras}. If $\cA$ is a simple vertex algebra and $G$ is a compact group of automorphisms of $\cA$, a theorem of Dong-Li-Mason asserts that $G$ and the invariant subalgebra $\cA^G$ behave like a dual reductive pair \DLM. As a module over $G$ and $\cA^G$, $\cA$ has a decomposition of the form \eqn\vaduality{\cA \cong \bigoplus_{\lambda\in I} W_{\lambda} \otimes V_{\lambda}.} Here $I$ indexes the set of irreducible, finite-dimensional representations $W_{\lambda}$ of $G$, and the $V_{\lambda}$'s are inequivalent, irreducible $\cA^G$-modules. In particular, $\cA^G$ is itself a simple vertex algebra, corresponding to the trivial $G$-module in \vaduality.

This is a beautiful and important result, but it tells us little about the {\it structure} of invariant vertex algebras of the form $\cA^G$. A range of examples of such invariant vertex algebras have been studied in both the physics and mathematics literature \BFH\EFH\KWY, and a general picture of their behavior has emerged. Often, $\cA$ is a current algebra, $bc$-system, $\beta\gamma$-system, or a tensor product of such algebras, and $\cA$ admits a $G$-invariant filtration for which $gr(\cA)$ is a commutative algebra with a derivation (i.e., an abelian vertex algebra). In this case, $gr(\cA^G) \cong (gr(\cA))^G$ as commutative algebras, and $gr(\cA^G)$ is viewed as a \lq\lq classical limit" of $\cA^G$. A key property is that $\cA^G$ tends to be {\it strongly finitely generated} as a vertex algebra. In other words, there exists a finite set of generators such that $\cA^G$ is spanned by the set of iterated Wick products of the generators and their derivatives. On the other hand, $gr(\cA^G)$ is usually {\it not} finitely generated as a vertex algebra, and a presentation of $gr(\cA^G)$ generally requires both infinitely many generators and infinitely many relations. The first example in \BFH, in which $\cA$ is a $\beta\gamma$-system and $G=SL_2$, nicely illustrates this phenomenon. Many interesting vertex algebras (such as various $\cW$-algebras) can be realized as invariant vertex algebras, but to the best of our knowledge there are no general theorems asserting that vertex algebras of the form $\cA^G$ are strongly finitely generated.

In this paper, we focus on invariant subalgebras of the $\beta\gamma$-system $\cS(V)$, or algebra of chiral differential operators, associated to the vector space $V = \C^n$. The standard action of $GL_n$ on $V$ induces an action of $GL_n$ on $\cS(V)$ by automorphisms, and we may consider $\cS(V)^G$ for any reductive subgroup $G\subset GL_n$. In the case $G=GL_n$, $\cS(V)^{GL_n}$ was shown by Kac-Radul to be isomorphic to the vertex algebra $\cW_{1+\infty,-n}$ with central charge $-n$ \KR. For any $c\in\C$, $\cW_{1+\infty,c}$ may be regarded as a module over the universal central extension $\hat{\cD}$ of the Lie algebra of differential operators on the circle. These algebras have been studied extensively in the physics literature, and they also play an important role in the theory of integrable systems. For $n=-1$, $\cW_{1+\infty,-1}$ is isomorphic to $\cW(\gg\gl_3)$ with central charge $-2$ \WI\WII, but less is known about the structure of $\cW_{1+\infty,-n}$ for $n>1$. It was conjectured in \BEHHH~and \WIV~that $\cW_{1+\infty,-n}$ should have a minimal, strong generating set consisting of $n^2+2n$ generators, and this conjecture was recently settled in \LII.

Our main result is that for any $G$ and $V$ as above, $\cS(V)^G$ is strongly finitely generated as a vertex algebra. The proof uses a combination of tools from classical invariant theory, together with the structure and representation theory of $\cW_{1+\infty,-n}$. As a vector space, $\cS(V)^G$ is isomorphic to the classical invariant ring $$R=(Sym \bigoplus_{k\geq 0} (V_k\oplus V^*_k))^{G},$$ where $V_k$ and $V^*_k$ are copies of $V$ and $V^*$, respectively. We view $\cS(V)^{G}$ as a {\it deformation} of $R$, in the sense that $\cS(V)^{G}$ admits a filtration for which the associated graded object $gr(\cS(V)^G)$ is isomorphic to $R$ as a commutative ring. In the terminology of Weyl, a {\it first fundamental theorem of invariant theory} for the pair $(G,V)$ is a set of generators for $R$. Even though $R$ is not finitely generated, only finitely many different \lq\lq types" of generators are necessary, and all others may be obtained from this set by polarization \We. Using this fact, together with the decomposition of $\cS(V)$ as a bimodule over $GL_n$ and $\cW_{1+\infty,-n}$ \KR, it is immediate that $\cS(V)^G$ is finitely generated as a vertex algebra. 

It takes more work to produce a strong finite generating set for $\cS(V)^G$. In Weyl's terminology, a {\it second fundamental theorem of invariant theory} for $(G,V)$ is a generating set for the ideal of relations among the generators of $R$. For $G=GL_n$, the strong finite generation of $\cS(V)^{GL_n}$ is essentially a consequence of the second fundamental theorem for the standard representation of $GL_n$, and our proof in \LII~depends on a detailed description of $R$ in this case. Remarkably, we will establish the strong finite generation of $\cS(V)^G$ for arbitrary $G$ without an explicit description of $R$. The only additional ingredient we need is a certain finiteness property possessed by any irreducible, highest-weight $\cW_{1+\infty,-n}$-submodule of $\cS(V)$. Using similar methods, we also show that vertex algebras of the form $\cE(V)^G$ and $(\cS(V)\otimes \cE(V))^G$ are strongly finitely generated, where $\cE(V)$ is the $bc$-system or semi-infinite exterior algebra associated to $V$, and $G$ is a reductive subgroup of $GL_n$. Finally, we establish the strong finite generation of a family of commutant subalgebras of $\cS(V)$ under Heisenberg algebra actions. This strengthens a result from \LI~which asserts that these vertex algebras are finitely generated.

It is well known that if a vertex algebra $\cV$ is strongly generated by a set $\{\alpha_i(z)|~i\in I\}$, the Zhu algebra of $\cV$ is generated by $\{a_i|~i\in I\}$, where $a_i$ is the image of $\alpha_i(z)$ under the Zhu map \Zh. Hence the Zhu algebras of all the invariant vertex algebras and commutant vertex algebras considered in this paper are finitely generated. Describing the Zhu algebras of these vertex algebras is the first step towards understanding their representation theory.

\subsec{Acknowledgment} I thank B. Bakalov for suggesting to me that invariant vertex algebras of the form $\cS(V)^G$ should be finitely generated.

\newsec{Vertex algebras}
In this section, we define vertex algebras, which have been discussed from various different points of view in the literature \Bor \FHL\FLM \Kac\Li \LZII. We will follow the formalism developed in \LZII~and partly in \Li. Let $V=V_0 \oplus V_1$ be a super vector space over $\C$, and let $z,w$ be formal variables. The space $QO(V)$ of {\it quantum operators} or {\it fields} on $V$ denotes the space of all linear maps $$V\ra V((z)):=\{\sum_{n\in\Z} v(n) z^{-n-1}|
v(n)\in V,\ v(n)=0\ for\ n>>0 \}.$$ Each element $a\in QO(V)$ can be
uniquely represented as a power series
$$a=a(z):=\sum_{n\in\Z}a(n)z^{-n-1}\in (End\ V)[[z,z^{-1}]].$$ We
refer to $a(n)$ as the $n$-th Fourier mode of $a(z)$. Each $a\in
QO(V)$ is assumed to be of the shape $a=a_0 +a_1$ where $a_i:V_j\ra V_{i+j}((z))$ for $i,j\in\Z/2\Z$, and we write $|a_i| = i$.

On $QO(V)$ there is a set of non-associative bilinear operations
$\circ_n$, indexed by $n\in\Z$, which we call the $n$-th circle
products. For homogeneous $a,b\in QO(V)$, they are defined by
$$
a(w)\circ_n b(w)=Res_z a(z)b(w)~\iota_{|z|>|w|}(z-w)^n-
(-1)^{|a||b|}Res_z b(w)a(z)~\iota_{|w|>|z|}(z-w)^n.
$$
Here $\iota_{|z|>|w|}f(z,w)\in\C[[z,z^{-1},w,w^{-1}]]$ denotes the
power series expansion of a rational function $f$ in the region
$|z|>|w|$. We usually omit the symbol $\iota_{|z|>|w|}$ and just
write $(z-w)^{-1}$ to mean the expansion in the region $|z|>|w|$,
and write $-(w-z)^{-1}$ to mean the expansion in $|w|>|z|$. It is
easy to check that $a(w)\circ_n b(w)$ above is a well-defined
element of $QO(V)$.

The non-negative circle products are connected through the {\it
operator product expansion} (OPE) formula.
For $a,b\in QO(V)$, we have \eqn\opeform{a(z)b(w)=\sum_{n\geq 0}a(w)\circ_n
b(w)~(z-w)^{-n-1}+:a(z)b(w):\ ,} which is often written as
$a(z)b(w)\sim\sum_{n\geq 0}a(w)\circ_n b(w)~(z-w)^{-n-1}$, where
$\sim$ means equal modulo the term $$
:a(z)b(w):\ =a(z)_-b(w)\ +\ (-1)^{|a||b|} b(w)a(z)_+\ .$$ Here
$a(z)_-=\sum_{n<0}a(n)z^{-n-1}$ and $a(z)_+=\sum_{n\geq
0}a(n)z^{-n-1}$. Note that $:a(z)b(z):$ is a well-defined element of
$QO(V)$. It is called the {\it Wick product} of $a$ and $b$, and it
coincides with $a(z)\circ_{-1}b(z)$. The other negative circle products
are related to this by
$$ n!~a(z)\circ_{-n-1}b(z)=\ :(\partial^n a(z))b(z):\ ,$$
where $\partial$ denotes the formal differentiation operator
$\frac{d}{dz}$. For $a_1(z),...,a_k(z)\in QO(V)$, the $k$-fold
iterated Wick product is defined to be
\eqn\iteratedwick{:a_1(z)a_2(z)\cdots a_k(z):\ =\ :a_1(z)b(z):~,}
where $b(z)=\ :a_2(z)\cdots a_k(z):\ $. We often omit the formal variable $z$ when no confusion will arise.

The set $QO(V)$ is a nonassociative algebra with the operations
$\circ_n$ and a unit $1$. We have $1\circ_n a=\delta_{n,-1}a$ for
all $n$, and $a\circ_n 1=\delta_{n,-1}a$ for $n\geq -1$. A linear subspace $\cA\subset QO(V)$ containing 1 which is closed under the circle products will be called a quantum operator algebra (QOA).
In particular $\cA$ is closed under $\partial$
since $\partial a=a\circ_{-2}1$. A subset $S=\{a_i|\ i\in I\}$ of $\cA$ is said to generate $\cA$ if any element $a\in\cA$ can be written as a linear
combination of nonassociative words in the letters $a_i$, $\circ_n$, for
$i\in I$ and $n\in\Z$. We say that $S$ {\it strongly generates} $\cA$ if any $a\in\cA$ can be written as a linear combination of words in the letters $a_i$, $\circ_n$ for $n<0$. Equivalently, $\cA$ is spanned by the collection $\{ :\partial^{k_1} a_{i_1}(z)\cdots \partial^{k_m} a_{i_m}(z):~|~ k_1,\dots,k_m \geq 0\}$. 

We say that $a,b\in QO(V)$ {\it quantum commute} if $(z-w)^N
[a(z),b(w)]=0$ for some $N\geq 0$. Here $[,]$ denotes the super bracket. This condition implies that $a\circ_n b = 0$ for $n\geq N$, so \opeform~becomes a finite sum. If $N$ can be chosen to be 0, we say that $a,b$ commute. A commutative quantum operator algebra (CQOA) is a QOA whose elements pairwise quantum commute. Finally, the notion of a CQOA is equivalent to the
notion of a vertex algebra. Every CQOA
$\cA$ is itself a faithful $\cA$-module, called the {\it left regular
module}. Define
$$\rho:\cA\rightarrow QO(\cA),\ \ \ \ a\mapsto\hat a,\ \ \ \ \hat
a(\zeta)b=\sum_{n\in\Z} (a\circ_n b)~\zeta^{-n-1}.$$ Then $\rho$ is an injective QOA homomorphism,
and the quadruple of structures $(\cA,\rho,1,\partial)$ is a vertex
algebra in the sense of \FLM. Conversely, if $(V,Y,{\bf 1},D)$ is
a vertex algebra, the collection $Y(V)\subset QO(V)$ is a
CQOA. {\it We will refer to a CQOA simply as a
vertex algebra throughout the rest of this paper}.

The following well-known identities measure the non-associativity and non-commutativity of the Wick product, and the failure of the positive circle products to be left and right derivations of the
Wick product. Let $a,b,c$ be vertex operators in some vertex algebra $\cA$, and let $n > 0$. Then
\eqn\vaidi{:(:ab:)c:-:abc:=\sum_{k\geq0}{1\over(k+1)!}\left(:(\partial^{k+1}a)(b\circ_k
c): +(-1)^{|a||b|}:(\partial^{k+1}b)(a\circ_k c):\right)}
\eqn\vaidii{:ab:-(-1)^{|a||b|}:ba:=\sum_{k\geq0}{(-1)^k\over(k+1)!}\partial^{k+1}(a\circ_kb),}
\eqn\vaidiii{a\circ_n(:bc:)-:(a\circ_nb)c:-(-1)^{|a||b|}:b(a\circ_nc):= \sum_{k=1}^n\left(\matrix{n\cr k}\right)(a\circ_{n-k}b)\circ_{k-1}c.}

\eqn\vaidiv{(:ab:)\circ_n
c=\sum_{k\geq0}{1\over k!}:(\partial^ka)(b\circ_{n+k}c):
+(-1)^{|a||b|}\sum_{k\geq0}b\circ_{n-k-1}(a\circ_k c) .}

\subsec{$\beta\gamma$-systems}

Let $V$ be a vector space of dimension $n$ over $\C$. The $\beta\gamma$-system $\cS(V)$, or algebra of chiral differential operators on $V$, was introduced by Friedan-Martinec-Shenker in \FMS. It is the unique even vertex algebra with generators $\beta^{x}(z)$, $\gamma^{x'}(z)$ for $x\in V$, $x'\in V^*$, which satisfy the OPE relations
$$\beta^x(z)\gamma^{x'}(w)\sim\langle x',x\rangle (z-w)^{-1},~~~~~~~\gamma^{x'}(z)\beta^x(w)\sim -\langle x',x\rangle (z-w)^{-1},$$
\eqn\betagamma{\beta^x(z)\beta^y(w)\sim 0,~~~~~~~\gamma^{x'}(z)\gamma^{y'}(w)\sim 0.} Here $\bra,\ket$ denotes the natural pairing between $V^*$ and $V$. We give $\cS(V)$ the conformal structure \eqn\virasoroelement{L(z) = \sum_{i=1}^n :\beta^{x_i}(z)\partial\gamma^{x'_i}(z):~,} under which $\beta^{x_i}(z)$ and $\gamma^{x'_i}(z)$ are primary of conformal weights $1$ and $0$, respectively. Here $\{x_1,\dots,x_n\}$ is a basis for $V$ and $\{x'_1,\dots,x'_n\}$ is the dual basis for $V^*$. 

The standard action $\rho:GL_n\ra Aut(V)$ induces an action $\hat{\rho}: GL_n\ra Aut(\cS(V))$ by vertex algebra automorphisms, defined on generators by \eqn\actionofg{\hat{\rho}(g)(\beta^x(z)) = \beta^{\rho(g)(x)}(z),~~~~ \hat{\rho}(g)( \gamma^{x'}(z)) = \gamma^{\rho^*(g)(x')}(z),~~~~ g\in G,~ x\in V,~ x'\in V^*.}

\newsec{Category $\cR$}
Let $\cR$ be the category of vertex algebras $\cA$ equipped with a $\Z_{\geq 0}$-filtration
\eqn\goodi{\cA_{(0)}\subset\cA_{(1)}\subset\cA_{(2)}\subset \cdots,\ \ \ \cA = \bigcup_{k\geq 0}
\cA_{(k)}} such that $\cA_{(0)} = \C$, and for all
$a\in \cA_{(k)}$, $b\in\cA_{(l)}$, we have
\eqn\goodii{a\circ_n b\in\cA_{(k+l)},\ \ \ for\
n<0,}
\eqn\goodiii{a\circ_n b\in\cA_{(k+l-1)},\ \ \ for\
n\geq 0.}
An element $a(z)\in\cA$ is said to have
degree $d$ if $d$ is the minimal integer for which
$a(z)\in\cA_{(d)}$. Morphisms in $\cR$ are vertex algebra homomorphisms which preserve the filtration.

Filtrations on vertex algebras satisfying \goodii-\goodiii~were introduced in \Lii~and are known as {\it good increasing filtrations}. If $\cA$ possesses such a filtration, the associated graded object $gr(\cA) = \bigoplus_{k>0}\cA_{(k)}/\cA_{(k-1)}$ is a
$\Z_{\geq 0}$-graded associative, supercommutative algebra with a
unit $1$ under a product induced by the Wick product on $\cA$. For $r\geq 1$, we denote by $\phi_r: \cA_{(r)} \ra \cA_{(r)}/\cA_{(r-1)}\subset gr(\cA)$ the natural projection. 
The operator $\partial = \frac{d}{dz}$ on $\cA$ induces a derivation $\partial$ of degree zero on $gr(\cA)$, and for each $a\in\cA_{(d)}$ and $n\geq 0$, the operator $a\circ_n$ on $\cA$
induces a derivation of degree $d-k$ on $gr(\cA)$. Here $$k  = sup \{ j\geq 1|~ \cA_{(r)}\circ_n \cA_{(s)}\subset \cA_{(r+s-j)}~\forall r,s,n\geq 0\},$$ as in \LL. The assignment $\cA\mapsto gr(\cA)$ is a functor from $\cR$ to the category of $\Z_{\geq 0}$-graded supercommutative rings with a differential $\partial$ of degree 0, which we will call $\partial$-rings. A $\partial$-ring is the same thing as an {\it abelian} vertex algebra, that is, a vertex algebra $\cV$ in which $[a(z),b(w)] = 0$ for all $a,b\in\cV$ \Bor. A $\partial$-ring $A$ is said to be generated by a subset $\{a_i|~i\in I\}$ if $\{\partial^k a_i|~i\in I, k\geq 0\}$ generates $A$ as a graded ring. The key feature of $\cR$ is the following reconstruction property \LL:

\lemma{Let $\cA$ be a vertex algebra in $\cR$ and let $\{a_i|~i\in I\}$ be a set of generators for $gr(\cA)$ as a $\partial$-ring, where $a_i$ is homogeneous of degree $d_i$. If $a_i(z)\in\cA_{(d_i)}$ are vertex operators such that $\phi_{d_i}(a_i(z)) = a_i$, then $\cA$ is strongly generated as a vertex algebra by $\{a_i(z)|~i\in I\}$.}\thmlab\recon

The main example we have in mind is $\cS(V)$, where we define $\cS(V)_{(r)}$ to be the linear span of the collection \eqn\goodsv{\{:\partial^{k_1} \beta^{x_1} \cdots \partial^{k_s}\beta^{x_s}\partial^{l_1} \gamma^{y'_1}\cdots \partial^{l_t}\gamma^{y'_t}:| ~~ x_i\in V,~y'_i\in V^*,~ k_i,l_i\geq 0,~ s+t \leq r\}.} Then $\cS(V)\cong gr(\cS(V))$ as linear spaces, and as a commutative algebra, we have
\eqn\structureofgrs{gr(\cS(V))\cong Sym \bigoplus_{k\geq 0} (V_k\oplus V^*_k),~~~~V_k = \{\beta^{x}_k |~ x\in V\},~~~~V^*_k = \{\gamma^{x'}_k |~ x'\in V^*\}.} In this notation, $\beta^{x}_k$ and $\gamma^{x'}_k$ are the images of $\partial^k \beta^{x}(z)$ and $\partial^k\gamma^{x'}(z)$ in $gr(\cS(V))$ under the projection $\phi_1: \cS(V)_{(1)}\ra \cS(V)_{(1)}/\cS(V)_{(0)}\subset gr(\cS(V))$. The action of $GL_n$ on $\cS(V)$ given by \actionofg~induces an action of $GL_n$ on $gr(\cS(V))$ by algebra automorphisms, and for all $k\geq 0$ we have isomorphisms of $GL_n$-modules $V_k\cong V$ and $V_k^*\cong V^*$. Finally, for any subgroup $G\subset GL_n$, $\cS(V)^G\cong gr(\cS(V)^G)$ as linear spaces, and \eqn\structureofgrsinv{gr(\cS(V)^G )\cong (gr(\cS(V))^G \cong (Sym \bigoplus_{k\geq 0} (V_k\oplus V^*_k))^G} as commutative algebras.

\newsec{The vertex algebra $\cW_{1+\infty,c}$}
Let $\cD$ be the Lie algebra of regular differential operators on $\C\setminus \{0\}$, with coordinate $t$. A standard basis for $\cD$ is $$J^l_k = -t^{l+k} (\partial_t)^l,~~~k\in \Z,~~~l\in \Z_{\geq 0},$$ where $\partial_t = \frac{d}{dt}$. $\cD$ has a 2-cocycle given by \eqn\cocycle{\Psi\bigg(f(t) (\partial_t)^m,  g(t) (\partial_t)^n\bigg) = \frac{m! n!}{(m+n+1)!} Res_{t=0} f^{(n+1)}(t) g^{(m)}(t) dt,} and a corresponding central extension $\hat{\cD} = \cD \oplus \C \kappa$, which was first studied by Kac-Peterson in \KP. $\hat{\cD}$ has a $\Z$-grading $\hat{\cD} = \bigoplus_{j\in\Z} \hat{\cD}_j$ by weight, given by
$$wt (J^l_k) = k,~~~~ wt( \kappa) = 0,$$ and a triangular decomposition $$\hat{\cD} = \hat{\cD}_+\oplus\hat{\cD}_0\oplus \hat{\cD}_-,$$ where $\hat{\cD}_{\pm} = \bigoplus_{j\in \pm \N} \hat{\cD}_j$ and $\hat{\cD}_0 = \cD_0\oplus \C\kappa.$
For a fixed $c\in\C$ and $\lambda\in \cD_0^*$, define the Verma module with central charge $c$ over $\hat{\cD}$ by
$$\cM_c(\hat{\cD},\lambda) = U(\hat{\cD})\otimes_{U(\hat{\cD}_0\oplus \hat{\cD}_+)} \C_{\lambda},$$ where $\C_{\lambda}$ is the one-dimensional $\hat{\cD}_0\oplus \hat{\cD}_+$-module on which $\kappa$ acts by multiplication by $c$ and $h\in\hat{\cD}_0$ acts by multiplication by $\lambda(h)$, and $\hat{\cD}_+$ acts by zero. There is a unique irreducible quotient of $\cM_c(\hat{\cD},\lambda)$ denoted by $V_c(\hat{\cD},\lambda)$.

Let $\cP$ be the parabolic subalgebra of $\cD$ consisting of differential operators which extend to all of $\C$, which has a basis $\{J^l_k|~l\geq 0,~l+k\geq 0\}$. The cocycle $\Psi$ vanishes on $\cP$, so $\cP$ may be regarded as a subalgebra of $\hat{\cD}$. Clearly $\hat{\cD}_0\oplus \hat{\cD}_+\subset \hat{\cP}$, where $\hat{\cP} = \cP\oplus \C\kappa$. The induced $\hat{\cD}$-module $$\cM_c=\cM_c(\hat{\cD},\hat{\cP}) = U(\hat{\cD})\otimes_{U(\hat{\cP})} \C_0$$ is then a quotient of $\cM_c(\hat{\cD},0)$, and is known as the {\it vacuum $\hat{\cD}$-module of central charge $c$}. $\cM_c$ has the structure of a vertex algebra which is freely generated by fields $$J^l(z) = \sum_{k\in\Z} J^l_k z^{-k-l-1},~~~~~~~l\geq 0$$ of weight $l+1$. The modes $J^l_k$ represent $\hat{\cD}$ on $\cM_c$, and we rewrite these fields in the form \eqn\defofjl{J^l(z) = \sum_{k\in\Z} J^l(k) z^{-k-1},~~~~~~~J^l(k) = J^l_{k-l}.}

An element $\omega\in \cM_c$ is called a {\it singular vector} if $J^l(k)(\omega) = 0$ for all $l\geq 0$ and $k>l$. The maximal proper $\hat{\cD}$-submodule $\cI_c$ is a vertex algebra ideal in $\cM_c$, and the unique irreducible quotient $\cM_c/\cI_c$ is denoted by $\cW_{1+\infty,c}$. In \LII, we denoted the projection $\cM_c\ra \cW_{1+\infty,c}$ by $\pi_{c}$, and we used the notation $j^l = \pi_c (J^l)$ in order to distinguish between $J^l\in \cM_{c}$ and its image in $\cW_{1+\infty,c}$. In this paper we only work with $\cW_{1+\infty,c}$, so no such distinction is necessary, and by abuse of notation we will denote the generators of $\cW_{1+\infty,c}$ by $J^l(z)$.

We are interested in the case of negative integral central charge. For $n\geq 1$, $\cW_{1+\infty,-n}$ has an important realization as a subalgebra of $\cS(V)$ for $V=\C^n$, which was introduced by Kac-Radul in \KR. It is given by
\eqn\bgrealization{J^l(z) \mapsto \sum_{i=1}^n :\gamma^{x'_i}(z) \partial^l \beta^{x_i}(z):~,} and the image of this embedding is precisely the invariant space $\cS(V)^{GL_n}$. Using this realization, together with Weyl's first and second fundamental theorems of invariant theory for the standard representation of $GL_n$, we showed in \LII~that $\cI_{-n}$ is generated by a singular vector of weight $(n+1)^2$. Moreover, this singular vector gives rise to a decoupling relation in $\cW_{1+\infty,-n}$ of the form \eqn\decoup{J^l = P(J^0,\dots,J^{l-1}),} for $l=n^2+2n$. Here $P$ is a normally ordered polynomial in the vertex operators $J^0,\dots, J^{l-1}$ and their derivatives. An easy consequence is that for all $r>l$, there exists a decoupling relation \eqn\higherdecoup{J^r = Q_r(J^0,\dots,J^{l-1}).} It follows that $\cW_{1+\infty,-n}$ is in fact {\it strongly} generated by $J^0,\dots, J^{l-1}$.

\newsec{Invariant subalgebras of $\beta\gamma$-systems}

Let $V = \C^n$, and let $G$ be a reductive subgroup of $GL_n$. Our goal is to describe the invariant subalgebra $\cS(V)^G$. Since $\cS(V)^{GL_n}\subset \cS(V)^G$, $\cS(V)^G$ is a module over $\cW_{1+\infty,-n}$, and this module structure will be an essential ingredient of our description.

\lemma{For any reductive $G\subset GL_n$, $\cS(V)^G$ is finitely generated as a vertex algebra.}\thmlab\sfg

\proof First we consider $\cS(V)^G$ from the point of view of classical invariant theory. Recall that $\cS(V)\cong gr(\cS(V))$ as linear spaces, and $$gr(\cS(V)^G )\cong (gr(\cS(V))^G \cong (Sym \bigoplus_{k\geq 0} (V_k\oplus V^*_k))^G = R$$ as commutative algebras. For any $p\geq 0$, there is an action of $GL_p$ on $\bigoplus_{k =0}^{p-1} (V_k\oplus V^*_k)$ which commutes with the action of $G$. The natural inclusions $GL_p\hookrightarrow GL_q$ for $p<q$ sending $$M \ra  \bigg[ \matrix{  M & 0 \cr 0 & I_{q-p} }\bigg]$$ induces an action of $GL_{\infty} = \lim_{p\ra \infty} GL_p$ on $\bigoplus_{k\geq 0} (V_k\oplus V^*_k)$. We obtain an action of $GL_{\infty}$ on $Sym \bigoplus_{k\geq 0} (V_k\oplus V^*_k)$ by algebra automorphisms, which commutes with the action of $G$. Hence $GL_{\infty}$ acts on $R$ as well. By a basic theorem of Weyl, $R$ is generated by the set of translates under $GL_{\infty}$ of any set of generators for $(Sym \bigoplus_{k = 0} ^{n-1}(V_k\oplus V^*_k))^G$ \We. Since $G$ is reductive, $(Sym \bigoplus_{k = 0} ^{n-1}(V_k\oplus V^*_k))^G$ is finitely generated. Hence there exists a finite set of homogeneous elements $\{f_1,\dots, f_k\}\subset R$ such that $\{ \sigma f_i|~ i=1,\dots,k,~ \sigma\in GL_{\infty}\}$ generates $R$. It follows from Lemma \recon~that the set of vertex operators $$\{(\sigma f_i)(z)\in \cS(V)^G|~i=1,\dots,k,~ \sigma\in GL_{\infty}\}$$ which correspond to $\sigma f_i$ under the linear isomorphism $\cS(V)^G\cong gr(\cS(V)^G) \cong R$, is a set of strong generators for $\cS(V)^G$.

Next, we recall the decomposition of $\cS(V)$ as a bimodule over $GL_n$ and $\cW_{1+\infty,-n}$, which appears in \KR. We have
\eqn\decompofsv{\cS(V) \cong \bigoplus_{\nu\in H} L(\nu)\otimes M^{\nu},} where $H$ indexes the irreducible, finite-dimensional representations $L(\nu)$ of $GL_n$, and $M^{\nu}$ is an irreducible, highest-weight $\cW_{1+\infty,-n}$-module. In particular, the $GL_n$-isotypic component of $\cS(V)$ of type $L(\nu)$ is isomorphic to $L(\nu)\otimes M^{\nu}$. Each $L(\nu)$ is a module over $G\subset GL_n$, and since $G$ is reductive, it has a decomposition $L(\nu) =\oplus_{\mu\in H^{\nu}} L(\nu)_{\mu}$. Here $\mu$ runs over a finite set $H^{\nu}$ of irreducible, finite-dimensional representations $L(\nu)_{\mu}$ of $G$, possibly with multiplicity. We thus obtain a refinement of \decompofsv:
\eqn\decompref{\cS(V) \cong \bigoplus_{\nu\in H} \bigoplus_{\mu\in H^{\nu}} L(\nu)_{\mu} \otimes M^{\nu}.}
Let $f_1(z),\dots,f_k(z)\in \cS(V)^G$ be the vertex operators corresponding to the polynomials $f_1, \dots,f_k$ under the linear isomorphism $\cS(V)^G\cong gr(\cS(V)^G) \cong R$. Clearly $f_1(z),\dots, f_k(z)$ must live in a finite direct sum
\eqn\newsummation{\bigoplus_{j=1}^r L(\nu_j)\otimes M^{\nu_j}} of the modules appearing in \decompofsv. By enlarging the collection $f_1(z),\dots,f_k(z)$ if necessary, we may assume without loss of generality that each $f_i(z)$ lives in a single representation of the form $L(\nu_j)\otimes M^{\nu_j}$. Moreover, we may assume that $f_i(z)$ lives in a trivial $G$-submodule $L(\nu_j)_{\mu_0} \otimes M^{\nu_j}$, where $\mu_0$ denotes the trivial, one-dimensional $G$-module. (In particular, $L(\nu_j)_{\mu_0}$ is one-dimensional). Since the actions of $GL_{\infty}$ and $GL_n$ on $\cS(V)$ commute, we may assume that $(\sigma f_i)(z)\in L(\nu_j)_{\mu_0}\otimes M^{\nu_j}$ for all $\sigma\in GL_{\infty}$. Since $\cS(V)^G$ is strongly generated by the set $\{ (\sigma f_i)(z)|~i=1,\dots,k,~ \sigma\in GL_{\infty}\}$, and each $M^{\nu_j}$ is an irreducible $\cW_{1+\infty,-n}$-module, $\cS(V)^G$ is generated as an algebra over $\cW_{1+\infty,-n}$ by $f_1(z),\dots,f_k(z)$. Finally, since $\cW_{1+\infty,-n}$ is itself a finitely generated vertex algebra, we conclude that $\cS(V)^G$ is finitely generated. $\Box$

We will refine this result to produce a {\it strong} finite generating set for $\cS(V)^G$. The fact that such a generating set exists in the case $G=GL_n$ is a formal consequence of Weyl's second fundamental theorem of invariant theory for the standard representation of $GL_n$ \LII. The decoupling relation \decoup~that gives rise to strong finite generation in this case is simply a deformation of the relation of minimal weight among the generators of $R=\big(Sym \bigoplus_{k\geq 0} (V_k\oplus V^*_k) \big)^{GL_n}$. We will establish the strong finite generation of $\cS(V)^G$ without any explicit knowledge of the second fundamental theorem for $(G,V)$. The only additional ingredient that we need is a certain finiteness property possessed by any irreducible, highest-weight $\cW_{1+\infty,-n}$-submodule of $\cS(V)$.

We begin with a basic observation about representations of associative algebras. Let $A$ be an associative $\C$-algebra (not necessarily unital), and let $W$ be a linear representation of $A$ (not necessarily finite-dimensional), via an algebra homomorphism $\rho: A\ra End(W)$. Regarding $A$ as a Lie algebra with commutator as bracket, let $\rho_{Lie}:A\ra End(W)$ denote the map $\rho$, regarded now as a Lie algebra homomorphism. There is then an induced algebra homomorphism $U(A)\ra End(W)$, where $U(A)$ denotes the universal enveloping algebra of $A$. Given elements $a,b\in A$, we denote the product in $U(A)$ by $a*b$ to distinguish it from $ab\in A$. Given a monomial $\mu = a_1* \cdots * a_r\in U(A)$, let $\tilde{\mu} = a_1\cdots a_r$ be the corresponding element of $A$. Let $U(A)_+$ denote the augmentation ideal (i. e., the ideal generated by $A$), regarded as an associative algebra with no unit. The map $U(A)_+ \ra A$ sending $\mu\mapsto \tilde{\mu}$ is then an algebra homomorphism which makes the diagram
 \eqn\commdiag{\matrix{
U(A)_+ & & \cr
\downarrow & \searrow \cr
A & \ra & End(W) \cr}}
commute. Let $T(W)$ denote the tensor algebra of $W$, whose $d$th graded component is the $d$-fold tensor product $W^{\otimes d}$. Clearly $\rho_{Lie}$ (but not $\rho$) can be extended to a Lie algebra homomorphism $\hat{\rho}_{Lie}: A\ra End(T(W))$, where $\hat{\rho}_{Lie}(a)$ acts by derivation on each $W^{\otimes d}$: $$\hat{\rho}_{Lie}(a)( w_1\otimes \cdots \otimes w_d) = \sum_{i=1}^d w_1\otimes \cdots \otimes \hat{\rho}_{Lie}(a)(w_i) \otimes \cdots \otimes w_d.$$ This extends to an algebra homomorphism $U(A)\ra End(T(W))$ which we also denote by $\hat{\rho}_{Lie}$, but there is no commutative diagram like \commdiag~because the map $A\ra End(T(W))$ is not a map of associative algebras. In particular, the restrictions of $\hat{\rho}_{Lie}(\mu)$ and $\hat{\rho}_{Lie}(\tilde{\mu})$ to $W^{\otimes d}$ are generally not the same for $d>1$.

\lemma{Given $\mu \in U(A)$ and $d\geq 1$, define a linear map $\Phi^d_{\mu} \in End(W^{\otimes d})$ by \eqn\mapmu{\Phi^d_{\mu}  = \hat{\rho}_{Lie}(\mu) \big|_{W^{\otimes d}} .} Let $E$ denote the subspace of $End(W^{\otimes d})$ spanned by $\{\Phi^d_{\mu}|~\mu\in U(A)\}$. Note that $E$ has a filtration $$E_1\subset E_2\subset \cdots,~~~~~~~E = \bigcup_{r\geq 1} E_r,$$ where $E_r$ is spanned by $\{\Phi^d_{\mu}|~ \mu \in U(A),~ deg(\mu) \leq r\}$. Then $E = E_d$.}\thmlab\first

\proof Given a monomial $\mu = a_1* \cdots * a_r \in U(A)$ of arbitrary degree $r>d$, we need to show that $\Phi^d_{\mu}$ can be expressed as a linear combination of elements of the form $\Phi^d_{\nu}$ where $\nu \in U(A)$ and $deg (\nu) \leq d$. Fix $p\leq d$, and let $Part^r_{p}$ denote the set of partitions $\phi$ of $\{1,\dots,r\}$ into $p$ disjoint, non-empty subsets $S^{\phi}_1,\dots, S^{\phi}_p$ whose union is $\{1,\dots,r\}$. Each subset $S^{\phi}_i$ is of the form $$S^{\phi}_i = \{i_1, \dots, i_{k_i} \},~~~~ i_1<\cdots < i_{k_i}.$$ For $i=1,\dots, p$, let $m_i\in U(A)$ be the corresponding monomial $m_i = a_{i_1} *\cdots *a_{i_{k_i}}$. Let $J = (j_1,\dots,j_p)$ be an (ordered) subset of $\{1,\dots,d\}$. Define a linear map $g_{\phi}\in End(W^{\otimes d})$ by \eqn\fphi{g_{\phi}(w_1\otimes \cdots \otimes w_d) =  \sum_J g_{\phi}^1(w_1) \otimes \cdots \otimes g_{\phi}^d (w_d),~~~~~~ g_{\phi}^k (w_k) =  \bigg\{ \matrix{ \hat{\rho}_{Lie}(m_i) (w_{j_i}) & k= j_i \cr & \cr w_k & k\neq j_i},}
where the sum runs over all (ordered) $p$-element subsets $J$ as above. Note that we could replace $m_i\in U(A)$ with $\tilde{m_i}\in A$ in \fphi, since $\hat{\rho}_{Lie}(m_i)(w_{j_i}) = \hat{\rho}_{Lie}(\tilde{m}_i)(w_{j_i})$. 

We claim that for each $\phi\in Part^r_p$, $g_{\phi} \in E_d$. We proceed by induction on $p$. The case $p=1$ is trivial because $g_{\phi} = \hat{\rho}_{Lie}(a)$ as derivations on $W^{\otimes d}$, where $a = a_1\cdots a_r$. Next, assume the result for all partitions $\psi\in Part^s_{q}$, for $q<p$ and $s\leq r$. Let $m_1,\dots,m_p\in U(A)$ be the monomials corresponding to $\phi$ as above, and define $m_{\phi} = \tilde{m}_1*\cdots * \tilde{m}_p \in U(A)$. By definition, $\Phi^d_{m_{\phi}} \in E_p\subset E_d$, and the leading term of $\Phi^d_{m_{\phi}}$ is $g_{\phi}$. The lower order terms are of the form $g_{\psi}$, where $\psi\in Part^p_q$ is a partition of $\{1,\dots, p\}$ into $q$ subsets, which each corresponds to a monomial in the variables $\tilde{m}_1,\dots,\tilde{m}_p$.
By induction, each of these terms lies in $E_q$, and since $g_{\phi} \equiv \Phi^d_{m_{\phi}}$ modulo $E_q$, the claim is proved.

Finally, using the derivation property of $A$ acting on $W^{\otimes d}$, one checks easily that \eqn\triplesum{\Phi^d_{\mu} = \sum_{p=1}^d \sum_{\phi\in Part^r_p} g_{\phi}.} Since each $g_{\phi}$ lies in $E_d$ by the above claim, this completes the proof of the lemma. $\Box$

A similar result clearly holds if we replace $T(W)$ with either $Sym(W)$ or $\bigwedge(W)$. We focus here on the symmetric case. First, the map $\rho_{Lie}: A\ra End(W)$ gives rise to a Lie algebra map $A\ra End(Sym(W))$, which extends to an algebra map $\hat{\rho}_{Lie}: U(A)\ra End(Sym(W))$. For each $d\geq 1$ and $\mu\in U(A)$, $\hat{\rho}_{Lie}(\mu)$ acts by derivation on the $d$th graded component $Sym^d(W)$, and we denote the restriction $\hat{\rho}_{Lie}(\mu) \big|_{Sym^d(W)}$ by $\Phi^d_{\mu}$. Let $E\subset End(Sym^d(W))$ be the vector space spanned by $\{\Phi^d_{\mu}|~\mu\in U(A)\}$, which has a filtration $E_1\subset E_2 \subset \cdots$, where $E_r$ is spanned by $\Phi^d_{\mu}$ for $\mu\in U(A)$ of degree at most $r$. For each partition $\phi\in Part^r_p$, \fphi~makes sense if the tensor product is replaced by the symmetric product, and \triplesum~holds. The same argument then shows that $E = E_d$.

\corollary{Let $f\in Sym^d(W)$, and let $M\subset Sym^d(W)$ be the cyclic $U(A)$-module generated by $f$. Then $\{\hat{\rho}_{Lie}(\mu)(f)|~ \mu\in U(A),~ deg(\mu)\leq d\}$ spans $M$.}\thmlab\firstcor

Recall the parabolic Lie subalgebra $\cP\subset \hat{\cD}$ with basis $\{J^l(k)|~k\geq 0\}$. We have a decomposition \eqn\decompofp{\cP = \cP_{-}\oplus \cD_0 \oplus \cP_+,~~~~~~~~~~~ \cP_{\pm} = \hat{\cD}_{\pm}\cap \cP.} In particular, $\cP_-$ has a basis $\{J^l(k)|~ 0\leq k<l\}$. Let $\cM$ be an irreducible, highest-weight $\cW_{1+\infty,-n}$-submodule of $\cS(V)$, with highest-weight vector $f(z)$, and let $\cM'$ be the $\cP$-submodule of $\cM$ generated by $f(z)$. By the Poincare-Birkhoff-Witt theorem, $\cM'$ is a quotient of $U(\cP)\otimes_{U(\cD_0\oplus \cP_+)} \C f(z)$, and in particular is a cyclic $\cP_-$-module with generator $f(z)$. Suppose that $f(z)$ has degree $d$, that is, $f(z)\in \cS(V)_{(d)} \setminus \cS(V)_{(d-1)}$. Since each element of $\cP$ preserves the filtration, and $\cM$ is irreducible, it is easy to see that the nonzero elements of $\cM'$ lie in $\cS(V)_{(d)} \setminus \cS(V)_{(d-1)}$. Therefore the projection $\cS(V)_{(d)}\ra \cS(V)_{(d)}/\cS(V)_{(d-1)} \subset gr(\cS(V))$ restricts to an isomorphism of $\cP$-modules
\eqn\isopmod{\cM'\cong gr(\cM')\subset gr(\cS(V)).}

\lemma{Let $\cM$ be an irreducible, highest-weight $\cW_{1+\infty,-n}$-submodule of $\cS(V)$ with highest-weight vector $f(z)$ of degree $d$. Let $\cM'$ be the corresponding $\cP$-module generated by $f(z)$. Then $\cM'$ is spanned by elements of the form $$\{J^{l_1}(k_1)\cdots J^{l_r}(k_r) f(z) |~ J^{l_i}(k_i)\in \cP_-,~ r\leq d\}.$$}\thmlab\second

\proof By \isopmod, $\cM'$ is isomorphic to the cyclic $\cP_-$-module $M = gr(\cM')$ generated by the image $f$ of $f(z)$ in $gr(\cS(V))$, which is homogeneous of degree $d$. The claim then follows from Corollary \firstcor, taking $A$ to be $\cP_-$ and $W$ to be the vector space with basis $\{\beta^{x_i}_k,\gamma^{x'_i}_k|~k\geq 0\}$. $\Box$

We need another fact about the structure of $\cS(V)$ as a module over $\cP$. For simplicity of notation, we take $n=1$, but the lemma we are going to prove holds for any $n$. In this case, $V= \C$ and $\cS(V)$ is generated by $\beta(z) = \beta^x(z)$ and $\gamma(z) = \gamma^{x'}(z)$. Let $W\subset gr(\cS(V))$ be the vector space with basis $\{\beta_k,\gamma_k|~ k\geq 0\}$, and for each $m\geq 0$, let $W_m$ be the subspace with basis $\{\beta_k,\gamma_k|~ 0\leq k\leq m\}$. Let $\phi:W\ra W$ be a linear map of weight $w\geq 1$, such that \eqn\arbmap{\phi(\gamma_i) = c_i \gamma_{i+w},~~~~~~~\phi(\beta_i) = d_i \beta_{i+w},} for constants $c_i,d_i\in \C$. For example, the restriction $J^{w+k}(k)\big|_{W}$ of any $J^{w+k}(k)\in \cP$, is such a map.

\lemma{Let $\phi$ be a linear map satisfying \arbmap. For any $m\geq 0$, the restriction $\phi \big|_{W_m}$ can be expressed uniquely as a linear combination of the operators $J^{w+k}(k)\big|_{W_m}$ for $0\leq k\leq 2m+1$.}\thmlab\third

\proof First, we need a basic calculation. For $w\geq 1$ and $k\geq 0$, we have
 \eqn\actionofp{J^{w+k}(k)(\beta_l) = \lambda^w_{k,l} (\beta_{l+w}),~~~~~~~J^{w+k}(k)(\gamma_l) = \mu^w_{k,l} (\gamma_{l+w}),$$ where $$\lambda^w_{k,l}  =  \bigg\{ \matrix{ - \frac{l!}{(l-k)!} & l-k \geq 0 \cr & \cr 0 & l-k <0},~~~~~~~\mu^w_{k,l}  =  (-1)^{w+k} \frac{(w+k+l)!}{(l+w)!}.}

Let $M^w$ be the block matrix  $$\bigg[ \matrix{  A^w & B^w \cr  C^w & D^w}\bigg],$$ where $A^w,B^w,C^w,D^w$ are $(m+1)\times(m+1)$ matrices with entries $$A_{ij} = \mu^w_{j,i},~~~~~B_{ij} = \mu^w_{m+1+j,i},~~~~~C_{ij} = \lambda^w_{j,i},~~~~~D_{ij} = \lambda^w_{m+1+j,i},~~~~0\leq i,j\leq m.$$ 
Let ${\bf c}$ be the column vector in $\C^{2m+2}$ whose transpose is given by $(c_0,\dots,c_m,d_0,\dots,d_m)$. Given an arbitrary linear combination $$\psi = t_0 J^w(0) + t_1 J^{w+1}(1) + \cdots + t_{2m+1} J^{w+2m+1}(2m+1)$$ of the operators $J^{w+k}(k)$ for $0\leq k\leq 2m+1$, let ${\bf t}$ be the column vector whose transpose is $(t_0,\dots, t_{2m+1})$. Note that $\phi \big|_{W_m} = \psi \big|_{W_m}$ precisely when $M^w  {\bf t} =  {\bf c}$, so in order to prove the claim, it suffices to show that $M^w$ is invertible. By \actionofp, $D^w$ is the zero matrix and $C^w$ is lower triangular with diagonal entries $C_{kk} = -k!$, so it suffices to show that $B^w$ is invertible. By \actionofp, the entries in each column of $B^w$ have the same sign, so $B^w$ is column-equivalent to the matrix \eqn\colequiv{ \left[\matrix{\frac{r!}{w!} & \frac{(r+1)!}{w!}  & \cdots & \frac{(r+m)!}{w!}  \cr & \cr  \frac{(r+1)!}{(w+1)!} &  \frac{(r+2)!}{(w+1)!} &  \cdots & \frac{(r+m+1)!}{(w+1)!} \cr & \cr  \vdots & \vdots & & \vdots  \cr & \cr  \frac{(r+m)!}{(w+m)!} & \frac{(r+m+1)!}{(w+m)!}  & \cdots & \frac{(r+2m)!}{(w+m)!}  
}\right],} where $r=w+m+1$. Clearly \colequiv~is row-equivalent to the matrix 
$$\left[ \matrix{1 & r+1  & (r+1)(r+2) & \cdot & \cdot & \cdot &  (r+1)(r+2)\cdots (r+m) \cr   1   & r+2  & (r+2)(r+3) & \cdot & \cdot &   \cdot & (r+2)(r+3)\cdots (r+m+1)  \cr  \vdots & \vdots & \vdots &  & & & \vdots  \cr  1  & r+m+1  & (r+m+1)(r+m+2) & \cdot & \cdot  &  \cdot & (r+m+1)(r+m+2)\cdots (r+2m) }\right],$$ which we denote by $T(r,m)$. We claim that $T(r,m)$ is invertible for any $r,m\geq 1$. For $m=1$, $T(r,1) = \bigg[ \matrix{  1 & r+1 \cr  1 & r+2 }\bigg]$, which is clearly invertible, so we may proceed by induction on $m$. For $m>1$, by subtracting the $i$th row from the $(i+1)$th row, for $i=0,\dots,m$, we see that $T(r,m)$ is row-equivalent to a matrix $\bigg[ \matrix{  1 & * \cr  0 & S }\bigg]$. Moreover, the $m\times m$ block $S$ is easily seen to be column-equivalent to $T(r+1,m-1)$. By our inductive assumption, $S$ is then invertible, so $T(r,m)$ is invertible as well. $\Box$

\remark{The same result holds for any $n\geq 1$. For $V = \C^n$, let $W\subset gr(\cS(V))$ be the vector space with basis $\{\beta^{x_i}_k,\gamma^{x'_i}_k|~i=1,\dots,n,~k\geq 0\}$, and let $W_m\subset W$ be the subspace with basis $\{\beta^{x_i}_k,\gamma^{x'_i}_k| ~ i=1,\dots,n,~0\leq k\leq m\}$. Let $\phi:W\ra W$ be a linear map of weight $w\geq1$ taking \eqn\actiononw{\gamma^{x'_i}_k\mapsto c_k \gamma^{x'_i}_{k+w},~~~~~ \beta^{x_i}_k\mapsto d_k \beta^{x_i}_{k+w},~~~~~i=1,\dots,n,} where the constants $c_k,d_k$ are independent of $i$. For example, each operator $\phi = J^{w+k}(k)\big|_W$ satisfies \actiononw. Then $\phi\big|_{W_m}$ can be expressed uniquely as a linear combination of $J^{w+k}(k)\big|_{W_m}$ for $0\leq k\leq 2m+1$. }\thmlab\remn

\lemma{Let $\cM$ be an irreducible, highest-weight $\cW_{1+\infty,-n}$-submodule of $\cS(V)$ with highest-weight vector $f(z)$ of degree $d$. Let $\cM'$ be the corresponding $\cP$-module generated by $f(z)$, and let $f$ be the image of $f(z)$ in $gr(\cS(V))$, which generates $M = gr(\cM')$ as a $\cP$-module. Fix $m$ so that $f\in Sym^d(W_m)$. Then $\cM'$ is spanned by $$\{J^{l_1}(k_1) \cdots J^{l_r}(k_r) f(z)|~J^{l_i}(k_i)\in \cP_-,~~ r\leq d,~~0\leq k_i \leq 2m+1\}.$$}\thmlab\fourth

\proof By \actionofp, we may work with $M = gr(\cM')$ rather than $\cM'$. For notational convenience, we do not distinguish between elements of $U(\cP_-)$ and their images in $End(Sym^d(W))$. As in the proof of Lemma \first, let $E$ denote the subspace of $End(Sym^d (W))$ spanned by $U(\cP_-)$, and let $E_r$ be the subspace spanned by elements of $U(\cP_-)$ of degree at most $r$. Let $\tilde{E}_r$ be the subspace of $E_r$ spanned by elements of $U(\cP_-)$ which only depend on $J^l(k)$ for $k\leq 2m+1$.

It is certainly not true that $E_d = \tilde{E}_d$ as subspaces of $End(Sym^d(W))$. However, it suffices to show that these spaces of endomorphisms coincide when restricted to $Sym^d(W_m)$. Since $E = E_d$, and hence is spanned by monomials $\mu = a_1*\cdots * a_r\in U(\cP_-)$ of degree $r\leq d$, we have 
\eqn\triplesumii{\mu = \sum_{p=1}^r  \sum_{\phi\in Part^r_p} g_{\phi},} where each partition $\phi\in Part^r_p$ corresponds to a set of monomials $m_1,\dots,m_p$, and $g_{\phi}$ is given by the symmetric version of \fphi, in which the tensor product is replaced by the symmetric product. For $p=r$, there is only one partition $\phi_0$ of $\{1,\dots,r\}$ into disjoint, non-empty subsets, and $g_{\phi_0}$ is defined on monomials $w_1\cdots w_d\in Sym^d(W)$ by \eqn\fphii{g_{\phi_0}(w_1\cdots w_d) =  \sum_J g_{\phi_0}^1(w_1)  \cdots  g_{\phi_0}^d (w_d),~~~~~~ g_{\phi_0}^k (w_k) =  \bigg\{ \matrix{ a_i (w_{j_i}) & k= j_i \cr & \cr w_k & k\neq j_i},} where the sum runs over all (ordered) $r$-element subsets $J \subset \{1,\dots,d\}$. By Lemma \third~and Remark \remn, the restriction of $a_i$ to $W_m$ coincides with a linear combination $S_i$ of the elements $J^l(k)\big|_{W_m}$ for $k\leq 2m+1$. Replace each of the factors $a_i (w_{j_i})$ appearing in \fphii~with $S_i (w_{j_i})$, and let $Q = \prod_{i=1}^r S_i$, which lies in $U(\cP_-)$, and depends only on $J^l(k)$ for $k\leq 2m+1$. Clearly the restriction of $Q$ to $Sym^d(W_m)$ agrees with the restriction of $\mu$ to $Sym^d(W_m)$, modulo terms lying in $E_{r-1}$. The lemma then follows by induction on $r$. $\Box$.

\remark{We may order the elements $J^l(k)\in \cP_-$ as follows: $J^{l_1}(k_1) > J^{l_2}(k_2)$ if $l_1>l_2$, or $l_1=l_2$ and $k_1<k_2$. Then Lemma \fourth~can be strengthened as follows: $\cM'$ is spanned by elements of the form $J^{l_1}(k_1)\cdots J^{l_r}(k_r) f(z)$ with \eqn\shapealpha{J^{l_i}(k_i)\in \cP_-,~~~~ r\leq d,~~~~0\leq k_i\leq 2m+1,~~~~ J^{l_1}(k_1)\geq \cdots \geq J^{l_r}(k_r).} The proof is a straightforward modification of the proof of Lemma \fourth, and is left to the reader.}\thmlab\secrem

In the next two lemmas, we use the notation $\cW_{1+\infty,-n}[k]$, $\cM[k]$, and $\cM'[k]$ to denote the homogeneous components of these spaces of conformal weight $k$. 

\lemma{Let $\cM$ be an irreducible, highest-weight $\cW_{1+\infty,-n}$-submodule of $\cS(V)$ with highest-weight vector $f(z)$. Define the Wick ideal of $\cM$ to be the subspace spanned by elements of the form $$:a(z) b(z):,~~~~~a(z)\in \bigoplus_{k>0} \cW_{1+\infty,-n}[k],~~~~~b(z)\in \cM.$$ Then any homogeneous element of $\cM$ of sufficiently high weight lies in the Wick ideal.}\thmlab\fifth

\proof It suffices to show that $\cM'[k]$ lies in the Wick ideal for $k>>0$, where $\cM'$ is the $\cP$-module generated by $f(z)$. As usual, let $d$ be the degree of $f(z)$, and fix $m$ so that $f \in Sym^d(W_m)$. By Remark \secrem, $\cM'$ is spanned by elements of the form $J^{l_1}(k_1)\cdots J^{l_r}(k_r) f(z)$ satisfying \shapealpha. Fix an element $\alpha(z)$ of this form of weight $K>>0$. Since each operator $J^{l_i}(k_i)$ has weight $l_i-k_i$, $k_i\leq 2m+1$, and $K>>0$, we may assume that $l_1>>n^2+2n$. Then \higherdecoup~allows us to express $J^{l_1}(z)$ as a normally ordered polynomial $Q_{l_1}(z)$ in the generators \eqn\genera{\partial^t J^l(z),~~~~~0\leq l<n^2+2n,~~~~~t\geq 0.} We claim that for any weight-homogeneous, normally ordered polynomial $Q(z)$ in the generators \genera~of sufficiently high weight, any element $c(z)\in \cM$, and any $k$ satisfying $0\leq k\leq 2m+1$, $Q(z)\circ_k c(z)$ lies in the Wick ideal. Specializing this to the case $Q(z) = Q_{l_1}(z)$, $c(z) = J^{l_2}(k_2)\cdots J^{l_r}(k_r) f(z)$, and $k=k_1$, proves the lemma.
 
We may assume without loss of generality that $Q(z)=:a(z)b(z):$ where $a(z) = \partial^t J^l(z)$ for some $0\leq l<n^2 +2n$. Then using \vaidiv, and supressing the formal variable $z$, we have \eqn\appvaid{Q\circ_k c = \big(:ab:\big) \circ_{k} c = \sum_{r\geq0}{1\over r!}:(\partial^r a)(b\circ_{k+r}c):
+\sum_{r\geq 0}b\circ_{k-r-1}(a\circ_r c) .} Suppose first that $b = \lambda 1$ for some constant $\lambda$. Then $Q= \lambda \partial^t J^l$, and since $wt(Q) >>0$, we have $t>>0$. Hence $Q\circ_k= \lambda(\partial^t J^l)\circ_k = 0$ as an operator (since this operator vanishes whenever $t>k$). So we may assume without loss of generality that $b$ is not a constant. 

We proceed by induction on $k$. For $k=0$, each term appearing in \appvaid~lies in the Wick ideal, so there is nothing to prove. For $k>0$, the only terms appearing in \appvaid~that need not lie in the Wick ideal a priori, are those of the form $\sum_{r=0}^{k-1} b\circ_{k-r-1}(a\circ_r c)$. However, each of these terms is weight-homogeneous, and the weight of $a\circ_r c = \partial^t J^l \circ_r c$ is bounded above by $wt(c) + n^2+2n$, since $\partial^t J^l \circ_r c=0$ for $t>r$. So we may still assume that $wt(b)>>0$. By our inductive assumption, all these terms then lie in the Wick ideal. $\Box$

Let $\cM$ be an irreducible, highest-weight $\cW_{1+\infty,-n}$-submodule of $\cS(V)$. Given a set $S$ of vertex operators in $\cM$, let $\cM_S\subset \cM$ denote the subspace spanned by elements of the form $:\omega_1(z)\cdots \omega_t(z) \alpha(z):$ with $\omega_j(z)\in \cW_{1+\infty,-n}$ and $\alpha(z)\in S$.

\lemma{Let $\cM$ be an irreducible, highest-weight $\cW_{1+\infty,-n}$-submodule of $\cS(V)$ with highest-weight vector $f(z)$. Then there exists a finite set $S$ of vertex operators of the form $$J^{l_1}(k_1) \cdots J^{l_r}(k_r) f(z),~~~~~~~~ 0\leq k_i < l_i<n^2+2n,$$ such that $\cM = \cM_S$.}\thmlab\sixth

\proof First, given a vertex operator of the form \eqn\eltofmp{J^{l_1}(k_1)\cdots J^{l_r}(k_r) f(z)\in \cM,} we can eliminate any operators of the form $J^{l_i}(k_i)$, $l_i\geq n^2+2n$ using the decoupling relations \decoup~and \higherdecoup~repeatedly, together with \vaidi-\vaidiv. Thus we can replace \eltofmp~with a linear combination of vertex operators of the form \eqn\replacement{:\partial^{a_1}J^{b_1}(z)\cdots\partial^{a_s}J^{b_s}(z) \big(J^{c_1}(m_1) \cdots J^{c_t}(m_t) f(z)\big) :,} where $b_i<n^2+2n$, and $0\leq m_j <c_j <n^2+2n$. 

Fix $K$ so that all homogeneous elements of $\cM$ of weight at least $K$ lie in the Wick ideal of $\cM$. Let $S$ denote the set of all elements of $\cM$ of the form \eqn\ess{\alpha(z) = J^{l_1}(k_1) \cdots J^{l_r}(k_r) f(z),~~~~~~ 0\leq k_i < l_i<n^2+2n,~~~~~~wt(\alpha(z))\leq K.} Clearly $S$ is a finite set. Given an element $\omega(z)\in\cM$ of higher weight, by applying Lemma \fifth~repeatedly, we can express $\omega(z)$ as a linear combination of elements of the form \eqn\desform{:a_1(z) \cdots a_r(z) b(z):,~~~~~a_i(z)\in \cW_{1+\infty,-n},~~~~~b(z)\in\bigoplus_{k\leq K} \cM[k].} Since each $b(z)$ has weight at most $K$, and is a linear combination of elements of the form \replacement, it follows that $b(z)\in \cM_S$. Then $\omega(z)\in \cM_S$ as well. $\Box$

Now we have assembled all the tools needed to prove our main result.

\theorem{Let $V=\C^n$ and let $G$ be a reductive subgroup of $GL_n$. Then $\cS(V)^G$ is strongly finitely generated as a vertex algebra.}\thmlab\strong

\proof By Lemma \sfg, we can find vertex operators $f_1(z),\dots, f_k(z)$ such that the corresponding polynomials $f_1,\dots, f_k\in gr(\cS(V))^G$, together with all $GL_{\infty}$ translates of $f_1,\dots, f_k$, generate the invariant ring $gr(\cS(V))^G$. As in the proof of Lemma \sfg, we may assume that each $f_i(z)$ lies in an irreducible, highest-weight $\cW_{1+\infty,-n}$-module $\cM_i$ of the form $L(\nu)_{\mu_0}\otimes M^{\nu}$, where $L(\nu)_{\mu_0}$ is a trivial, one-dimensional $G$-module. Furthermore, we may assume without loss of generality that $f_1(z),\dots, f_k(z)$ are highest-weight vectors for the action of $\cW_{1+\infty,-n}$. (Otherwise, we can replace these with the highest-weight vectors in the corresponding modules). 

For each $\cM_i$, we can choose a finite set $S_i$ of vertex operators of the form $$J^{l_1}(k_1)\cdots J^{l_r}(k_r) f_i(z),~~~~~~~~0\leq k_t <l_t< n^2+2n,$$ such that $\cM_i = (\cM_i)_{S_i}$. Define $$ S=\{J^0(z),\dots, J^{n^2+2n-1}(z) \} \cup \big(\bigcup_{i=1}^k S_i \big).$$ Since $\{J^0(z),\dots, J^{n^2+2n-1}(z)\}$ strongly generates $\cW_{1+\infty,-n}$, and the set $\bigcup_{i=1}^k \cM_i$ strongly generates $\cS(V)^G$, it is immediate that $S$ is a strong, finite generating set for $\cS(V)^G$. $\Box$

\newsec{Invariant subalgebras of $bc$-systems and $bc\beta\gamma$-systems}
Our methods easily extend to the study of invariant subalgebras of $bc$-systems and $bc\beta\gamma$-systems. Given a finite-dimensional vector space $V$, the $bc$-system $\cE(V)$, or semi-infinite exterior algebra, was introduced by Friedan-Martinec-Shenker in \FMS. It is the unique vertex algebra with odd generators $b^{x}(z)$, $c^{x'}(z)$ for $x\in V$, $x'\in V^*$, which satisfy the OPE relations
$$b^x(z)c^{x'}(w)\sim\langle x',x\rangle (z-w)^{-1},~~~~~~~c^{x'}(z) b^x(w)\sim \langle x',x\rangle (z-w)^{-1},$$
\eqn\bcsystem{b^x(z) b^y(w)\sim 0,~~~~~~~c^{x'}(z) c^{y'}(w)\sim 0.} The $bc\beta\gamma$ system on $V$ is defined to be $\cE(V)\otimes \cS(V)$. 

As shown in \FKRW, for $n\geq 1$, $\cW_{1+\infty,n}$ has a free field realization as the invariant subalgebra $\cE(V)^{GL_n}$ for $V=\C^n$. It is given by
\eqn\bcrealization{J^l(z) \mapsto \sum_{i=1}^n : c^{x'_i}(z) \partial^l b^{x_i}(z):~.} As a bimodule over $GL_n$ and $\cW_{1+\infty,n}$, $\cE(V)$ has a decomposition $$\cE(V) \cong \bigoplus_{\nu\in H} L(\nu)\otimes N^{\nu}$$ of the form \decompofsv, where $L(\nu)$ is an irreducible, finite-dimensional $GL_n$-module and $N^{\nu}$ is an irreducible, highest-weight $\cW_{1+\infty,n}$-module \KR. Hence the $bc\beta\gamma$-system $\cE(V)\otimes \cS(V)$ has a decomposition $$\cE(V)\otimes \cS(V)\cong \bigoplus_{\nu\in H,~ \mu\in H} L(\nu)\otimes L(\mu) \otimes M^{\nu}\otimes N^{\mu},$$ where $L(\nu)$ and $L(\mu)$ are irreducible, finite-dimensional $GL_n$-modules, and $M^{\nu}$ and $N^{\nu}$ are irreducible, highest-weight modules over $\cW_{1+\infty,-n}$ and $\cW_{1+\infty,n}$, respectively.

The same argument as the proof of Lemma \sfg~shows that for any reductive $G\subset GL_n$, both $\cE(V)^G$ and $(\cE(V)\otimes \cS(V))^G$ are finitely generated vertex algebras. Finally, the analogue of Lemma \sixth~holds for each irreducible, highest-weight $\cW_{1+\infty,n}$-submodule $\cM$ of $\cE(V)$ with highest-weight vector $f(z)$. Given a subset $S\subset \cM$, we define $\cM_S\subset \cM$ to be the subspace spanned by the elements $$:\omega_1(z)\cdots \omega_r(z) \alpha(z):~,~~~~\omega_j(z)\in \cW_{1+\infty,n},~~~~\alpha(z)\in S.$$ Then there is a finite set $S$ of vertex operators of the form $$ J^{l_1}(k_1)\cdots J^{l_r}(k_r)f(z),~~~~~ 0\leq k_i <l_i <n,$$ such that $\cM = \cM_S$. Similarly, for any irreducible $\cW_{1+\infty,n}\otimes \cW_{1+\infty,-n}$-submodule $\cM$ of $\cE(V)\otimes \cS(V)$ with highest-weight vector $f(z)$, and any subset $S\subset \cM$, define $\cM_S$ to be the subspace spanned by the elements $$:\omega_1(z)\cdots \omega_r(z) \nu_1(z)\cdots\nu_s(z) \alpha(z):$$ with $\omega_i(z)\in \cW_{1+\infty,n}$, $\nu_j(z)\in \cW_{1+\infty,-n}$, and $\alpha(z)\in S$. Then there is a finite set $S$ of vertex operators of the form
$$J^{l_1}(k_1)\cdots J^{l_r}(k_r) \tilde{J}^{d_1}(e_1)\cdots \tilde{J}^{d_s}(e_s) f(z),~~~~~ 0\leq k_i <l_i <n,~~~~~0\leq e_i <d_i <n^2+2n,$$ with $J^{l_i}\in \cW_{1+\infty,n}$ and $\tilde{J}^{d_j}\in\cW_{1+\infty,-n}$, such that $\cM = \cM_S$. An immediate consequence, whose proof is the same as the proof of Theorem \strong, is
\theorem{Let $V=\C^n$ and let $G$ be a reductive subgroup of $GL_n$. Then $\cE(V)^G$ and $(\cE(V)\otimes \cS(V))^G$ are strongly finitely generated as vertex algebras.}

We remark that $(\cE(V)\otimes \cS(V))^G$ has a natural nonlinear generalization. Let $X$ be a nonsingular algebraic variety over $\C$, equipped with an algebraic action of a reductive group $G$. In \MSV, Malikov-Schectman-Vaintrob introduced a sheaf of vertex algebras $\Omega^{ch}_X$ on $X$ known as the chiral de Rham sheaf. The algebra of global sections $\Omega^{ch}(X)$ admits an action of $G$ by automorphisms, and it is natural to study the invariant subalgebra $(\Omega^{ch}(X))^G$. In the case $X=V$, $\Omega^{ch}(V)$ is precisely $\cE(V)\otimes\cS(V)$, so if $G$ acts linearly on $V$, $(\Omega^{ch}(V))^G$ is strongly finitely generated. For a general $X$, the action of $\cW_{1+\infty,n}\otimes \cW_{1+\infty,-n}$ will not be globally defined on $\Omega^{ch}(X)$, so a new approach is needed to determine the structure of $(\Omega^{ch}(X))^G$.

Finally, we point out that for $V=\C^n$, both $\cE(V)$ and $\cS(V)$ admit additional automorphisms beyond those arising from the action of $GL_n$ on $V$ as above. For example, $SO_{2n}$ acts naturally on $\cE(V)$, and the decomposition of $\cE(V)$ as a bimodule over $SO_{2n}$ and $\cE(V)^{SO_{2n}}$ was described explicitly in \WIII. Moreover, $\cE(V)^{SO_{2n}}$ is closely related to the classical $\cW$-algebra $\cW\cD_n$ with central charge $n$ (see Theorem 14.2 of \KWY), and in particular is finitely generated. Likewise, $Sp_{2n}$ acts naturally on $\cS(V)$, and the decomposition of $\cS(V)$ as a bimodule over $Sp_{2n}$ and $\cS(V)^{Sp_{2n}}$ was also worked out in \WIII. We expect that the strong finite generation of $\cE(V)^{SO_{2n}}$ and $\cS(V)^{Sp_{2n}}$ can be established using the approach of \LII. (For example, the strong finite generation of $\cS(V)^{Sp_{2n}}$ should be a formal consequence of the second fundamental theorem of invariant theory for the standard representation of $Sp_{2n}$). Finally, we expect that the strong finite generation of $\cE(V)^G$ and $\cS(V)^{G'}$ for reductive subgroups $G\subset SO_{2n}$ and $G'\subset Sp_{2n}$, can be established using the method of this paper. This method can be summarized as follows: given an invariant vertex algebra $\cA^G$, find a \lq\lq big" subalgebra $\cB\subset \cA^G$ such that $\cB$ is strongly finitely generated, $\cA^G$ is completely reducible as a $\cB$-module, and $\cA^G$ is finitely generated as an algebra over $\cB$. If the irreducible $\cB$-submodules of $\cA^G$ have a finiteness property analogous to the property of the $\cW_{1+\infty,-n}$-submodules of $\cS(V)$ given by Lemma \sixth, $\cA^G$ will be strongly finitely generated.

\newsec{Torus actions and commutant subalgebras of $\cS(V)$}
Let $G$ be an algebraic torus of dimension $m$ acting faithfully and diagonally on $V = \C^n$. There is an induced representation $\rho:\gg \ra End(V)$, where $\gg$ is the (abelian) Lie algebra $\C^m$. In the notation of \LI, this induces a vertex algebra homomorphism $\hat{\tau}: \cO(\gg,B)\ra \cS(V)$, where $\cO(\gg,B)$ is the current algebra associated to $\gg$ equipped with the bilinear form $B(\xi,\eta) = -Tr(\rho(\xi)\rho(\eta))$. Since $\gg$ is abelian, $\cO(\gg,B)$ is just the tensor product of $m$ copies of the Heisenberg vertex algebra. In \LI, we studied the commutant subalgebra $Com(\hat{\tau}(\cO(\gg,B)),\cS(V))$, which is just the invariant space $\cS(V)^{\gg[t]}$. Philosophically, this problem is similar to studying invariant subalgebras of $\cS(V)$ under reductive group actions (in contrast to $\cS(V)^{\gg[t]}$ for nonabelian $\gg$), since the Heisenberg algebra acts semisimply. 

In \LI, we showed that $\cS(V)^{\gg[t]}$ is a finitely generated vertex algebra. First, $\cS(V)^{\gg[t]}$ contains a subalgebra $\cB' = \Phi \otimes \cW$, where $\Phi$ is the tensor product of $n-m$ copies of the Heisenberg algebra, and $\cW$ is the tensor product of $n$ copies of the Zamolodchikov $\cW_3$ algebra with central charge $-2$. Recall that $\cW_{3,-2}$ has generators $L$ and $W$ of weights $2$ and $3$, respectively. We denote the generators of $\Phi$ by $\phi^1,\dots,\phi^{n-m}$, and we denote the generators of $\cW$ by $L^1,W^1,\dots,L^n,W^n$. We have a direct sum decomposition \eqn\decompcomm{ \cS(V)^{\gg[t]} = \bigoplus_{l\in \cL} \cM'_l,} where $\cL$ is a certain lattice determined by the group action (denoted by $A^{\perp}\cap \Z^n$ in \LI), and $\cM'_l$ is the irreducible, cyclic $\cB'$-module with generator $\omega_l(z)$. If we choose a basis $l^1,\dots,l^r$ for the lattice $\cL$, the corresponding vertex operators $$\{\omega_{l^i}(z),\omega_{-l^i}(z)|~i=1,\dots,r\},$$ together with the generators of $\cB'$, are a finite generating set for $\cS(V)^{\gg[t]}$. However, this set is not generally a {\it strong} finite generating set.

\theorem{For any action of a torus $G$ on $V$ as above, $\cS(V)^{\gg[t]}$ is strongly finitely generated as a vertex algebra.}

\proof The main idea is that the $\cB'$-modules $\cM'_l$ appearing in the decomposition \decompcomm~of $\cS(V)^{\gg[t]}$ have a similar finiteness property to the one given by Lemma \sixth. For each basis element $l^i\in \cL$, there is a finite set $S_i$ of vertex operators in $\cM'_{l^i}$ of the form $$\big(L^{1}(0) \big)^{r_1} \big(W^{1}(0) \big)^{s_1} \big(W^1(1)\big)^{t_1} \cdots \big(L^{n}(0) \big)^{r_n} \big(W^{n}(0) \big)^{s_n} \big(W^n(1)\big)^{t_n} \omega_{l^i}(z)$$ such that $(\cM'_{l^i})_{S_i} = \cM'_{l^i}$. In this notation, $(\cM'_{l^i})_{S_i}$ denotes the subspace of $\cM'_{l^i}$ spanned by $$\{:\omega_1(z)\cdots\omega_t(z) \alpha(z):|~ \omega_j(z) \in \cB',~ \alpha(z)\in S_i\}.$$ The argument is similar to the proof of Lemma \sixth, and can in fact be obtained directly from Lemma \sixth~in the case $n=1$ by using the isomorphism $\cW_{1+\infty,-1}\cong \cH\otimes \cW_{3,-2}$ due to Wang \WI\WII. Similarly, there is a finite set $T_i$ of vertex operators in $\cM'_{-l^i}$ such that  $(\cM'_{-l^i})_{T_i} = \cM'_{-l^i}$. It is immediate that $$S=\{ \phi^i(z), L^j(z),W^j(z)|~i=1,\dots,n-m,~j=1,\dots,n\} \cup \bigcup_{k=1}^r (S_k \cup T_k),$$ is a strong finite generating set for $\cS(V)^{\gg[t]}$. $\Box$

\footatend\vfill\supereject\immediate\closeout\rfile\writestoppt
\baselineskip=14pt\centerline{{\bf References}}\bigskip{\frenchspacing%
\parindent=20pt\escapechar=` \input refs.tmp\vfill\eject}\nonfrenchspacing


\end

%% file: abbrev.tex

%
%
%
\def\unredoffs{} \def\redoffs{\voffset=-.31truein\hoffset=-.59truein}
\def\speclscape{\special{ps: landscape}}
%
%
%
%
\newbox\leftpage \newdimen\fullhsize \newdimen\hstitle \newdimen\hsbody
\tolerance=1000\hfuzz=2pt
\catcode`\@=11 
%
\ifx\answ\bigans\message{(This will come out unreduced.}
\magnification=1200\unredoffs\baselineskip=16pt plus 2pt minus 1pt
\hsbody=\hsize \hstitle=\hsize 
\else\message{(This will be reduced.} \let\l@r=L
\magnification=1000\baselineskip=16pt plus 2pt minus 1pt \vsize=7truein
\redoffs \hstitle=8truein\hsbody=4.75truein\fullhsize=10truein\hsize=\hsbody
\output={\ifnum\pageno=0 
  \shipout\vbox{\speclscape{\hsize\fullhsize\makeheadline}
    \hbox to \fullhsize{\hfill\pagebody\hfill}}\advancepageno
  \else
  \almostshipout{\leftline{\vbox{\pagebody\makefootline}}}\advancepageno
  \fi}
\def\almostshipout#1{\if L\l@r \count1=1 \message{[\the\count0.\the\count1]}
      \global\setbox\leftpage=#1 \global\let\l@r=R
 \else \count1=2
  \shipout\vbox{\speclscape{\hsize\fullhsize\makeheadline}
      \hbox to\fullhsize{\box\leftpage\hfil#1}}  \global\let\l@r=L\fi}
\fi
%
\newcount\yearltd\yearltd=\year\advance\yearltd by -1900

%

%
%

\def\draftmode{\message{ DRAFTMODE }\def\draftdate{{\rm preliminary draft:
\number\month/\number\day/\number\yearltd\ \ \hourmin}}%
\headline={\hfil\draftdate}\writelabels\baselineskip=20pt plus 2pt minus 2pt
 {\count255=\time\divide\count255 by 60 \xdef\hourmin{\number\count255}
  \multiply\count255 by-60\advance\count255 by\time
  \xdef\hourmin{\hourmin:\ifnum\count255<10 0\fi\the\count255}}}
\def\nolabels{\def\wrlabeL##1{}\def\eqlabeL##1{}\def\reflabeL##1{}}
\def\writelabels{\def\wrlabeL##1{\leavevmode\vadjust{\rlap{\smash%
{\line{{\escapechar=` \hfill\rlap{\sevenrm\hskip.03in\string##1}}}}}}}%
\def\eqlabeL##1{{\escapechar-1\rlap{\sevenrm\hskip.05in\string##1}}}%
\def\reflabeL##1{\noexpand\llap{\noexpand\sevenrm\string\string\string##1}}}
\nolabels
%
\global\newcount\secno \global\secno=0
\global\newcount\meqno \global\meqno=1
\def\newsec#1{\global\advance\secno by1\message{(\the\secno. #1)}
\global\subsecno=0\eqnres@t\noindent{\bf\the\secno. #1}
\writetoca{{\secsym} {#1}}\par\nobreak\medskip\nobreak}
\def\eqnres@t{\xdef\secsym{\the\secno.}\global\meqno=1\bigbreak\bigskip}
\def\sequentialequations{\def\eqnres@t{\bigbreak}}\xdef\secsym{}
\global\newcount\subsecno \global\subsecno=0
\def\subsec#1{\global\advance\subsecno by1\message{(\secsym\the\subsecno. #1)}
\ifnum\lastpenalty>9000\else\bigbreak\fi
\noindent{\it\secsym\the\subsecno. #1}\writetoca{\string\quad
{\secsym\the\subsecno.} {#1}}\par\nobreak\medskip\nobreak}
\def\appendix#1#2{\global\meqno=1\global\subsecno=0\xdef\secsym{\hbox{#1.}}
\bigbreak\bigskip\noindent{\bf Appendix #1. #2}\message{(#1. #2)}
\writetoca{Appendix {#1.} {#2}}\par\nobreak\medskip\nobreak}
%
%
\def\eqnn#1{\xdef #1{(\secsym\the\meqno)}\writedef{#1\leftbracket#1}%
\global\advance\meqno by1\wrlabeL#1}
\def\eqna#1{\xdef #1##1{\hbox{$(\secsym\the\meqno##1)$}}
\writedef{#1\numbersign1\leftbracket#1{\numbersign1}}%
\global\advance\meqno by1\wrlabeL{#1$\{\}$}}
\def\eqn#1#2{\xdef #1{(\secsym\the\meqno)}\writedef{#1\leftbracket#1}%
\global\advance\meqno by1$$#2\eqno#1\eqlabeL#1$$}
%
\newskip\footskip\footskip14pt plus 1pt minus 1pt 
\def\footnotefont{\ninepoint}\def\f@t#1{\footnotefont #1\@foot}
\def\f@@t{\baselineskip\footskip\bgroup\footnotefont\aftergroup\@foot\let\next}
\setbox\strutbox=\hbox{\vrule height9.5pt depth4.5pt width0pt}
\global\newcount\ftno \global\ftno=0
\def\foot{\global\advance\ftno by1\footnote{$^{\the\ftno}$}}
%
\newwrite\ftfile
\def\footend{\def\foot{\global\advance\ftno by1\chardef\wfile=\ftfile
$^{\the\ftno}$\ifnum\ftno=1\immediate\openout\ftfile=foots.tmp\fi%
\immediate\write\ftfile{\noexpand\smallskip%
\noexpand\item{f\the\ftno:\ }\pctsign}\findarg}%
\def\footatend{\vfill\eject\immediate\closeout\ftfile{\parindent=20pt
\centerline{\bf Footnotes}\nobreak\bigskip\input foots.tmp }}}
\def\footatend{}
%
%
\global\newcount\refno \global\refno=1
\newwrite\rfile
%
\def\ref{\nref}
\def\nref#1{\xdef#1{[\the\refno]}\writedef{#1\leftbracket#1}%
\ifnum\refno=1\immediate\openout\rfile=refs.tmp\fi
\global\advance\refno by1\chardef\wfile=\rfile\immediate
\write\rfile{\noexpand\item{#1\ }\reflabeL{#1\hskip.31in}\pctsign}\findarg}
\def\findarg#1#{\begingroup\obeylines\newlinechar=`\^^M\pass@rg}
{\obeylines\gdef\pass@rg#1{\writ@line\relax #1^^M\hbox{}^^M}%
\gdef\writ@line#1^^M{\expandafter\toks0\expandafter{\striprel@x #1}%
\edef\next{\the\toks0}\ifx\next\em@rk\let\next=\endgroup\else\ifx\next\empty%
\else\immediate\write\wfile{\the\toks0}\fi\let\next=\writ@line\fi\next\relax}}
\def\striprel@x#1{} \def\em@rk{\hbox{}}
\def\lref{\begingroup\obeylines\lr@f}
\def\lr@f#1#2{\gdef#1{\ref#1{#2}}\endgroup\unskip}

\def\addref#1{\immediate\write\rfile{\noexpand\item{}#1}} 
\def\startrefs#1{\immediate\openout\rfile=refs.tmp\refno=#1}
\def\refs#1{\count255=1[\r@fs #1{\hbox{}}]}
\def\r@fs#1{\ifx\und@fined#1\message{reflabel \string#1 is undefined.}%
\nref#1{need to supply reference \string#1.}\fi%
\vphantom{\hphantom{#1}}\edef\next{#1}\ifx\next\em@rk\def\next{}%
\else\ifx\next#1\ifodd\count255\relax\xref#1\count255=0\fi%
\else#1\count255=1\fi\let\next=\r@fs\fi\next}
%

%
\newwrite\ffile\global\newcount\figno \global\figno=1
\def\fig{fig.~\the\figno\nfig}
\def\nfig#1{\xdef#1{fig.~\the\figno}%
\writedef{#1\leftbracket fig.\noexpand~\the\figno}%
\ifnum\figno=1\immediate\openout\ffile=figs.tmp\fi\chardef\wfile=\ffile%
\immediate\write\ffile{\noexpand\medskip\noexpand\item{Fig.\ \the\figno. }
\reflabeL{#1\hskip.55in}\pctsign}\global\advance\figno by1\findarg}
\def\xfig{\expandafter\xf@g}\def\xf@g fig.\penalty\@M\ {}
\def\figs#1{figs.~\f@gs #1{\hbox{}}}
\def\f@gs#1{\edef\next{#1}\ifx\next\em@rk\def\next{}\else
\ifx\next#1\xfig #1\else#1\fi\let\next=\f@gs\fi\next}
\newwrite\lfile
{\escapechar-1\xdef\pctsign{\string\%}\xdef\leftbracket{\string\{}
\xdef\rightbracket{\string\}}\xdef\numbersign{\string\#}}

\def\writestop{\def\writestoppt{\immediate\write\lfile{\string\pageno%
\the\pageno\string\startrefs\leftbracket\the\refno\rightbracket%
\string\def\string\secsym\leftbracket\secsym\rightbracket%
\string\secno\the\secno\string\meqno\the\meqno}\immediate\closeout\lfile}}
\def\writestoppt{}\def\writedef#1{}
\def\seclab#1{\xdef #1{\the\secno}\writedef{#1\leftbracket#1}\wrlabeL{#1=#1}}
\def\subseclab#1{\xdef #1{\secsym\the\subsecno}%
\writedef{#1\leftbracket#1}\wrlabeL{#1=#1}}
\newwrite\tfile \def\writetoca#1{}
\def\leaderfill{\leaders\hbox to 1em{\hss.\hss}\hfill}
\def\writetoc{\immediate\openout\tfile=toc.tmp
   \def\writetoca##1{{\edef\next{\write\tfile{\noindent ##1
   \string\leaderfill {\noexpand\number\pageno} \par}}\next}}}
%
%
%

\catcode`\@=12 
%
\edef\tfontsize{\ifx\answ\bigans scaled\magstep3\else scaled\magstep4\fi}
\font\titlerm=cmr10 \tfontsize \font\titlerms=cmr7 \tfontsize
\font\titlermss=cmr5 \tfontsize \font\titlei=cmmi10 \tfontsize
\font\titleis=cmmi7 \tfontsize \font\titleiss=cmmi5 \tfontsize
\font\titlesy=cmsy10 \tfontsize \font\titlesys=cmsy7 \tfontsize
\font\titlesyss=cmsy5 \tfontsize \font\titleit=cmti10 \tfontsize
\skewchar\titlei='177 \skewchar\titleis='177 \skewchar\titleiss='177
\skewchar\titlesy='60 \skewchar\titlesys='60 \skewchar\titlesyss='60
\def\titlefont{\def\rm{\fam0\titlerm}
\textfont0=\titlerm \scriptfont0=\titlerms \scriptscriptfont0=\titlermss
\textfont1=\titlei \scriptfont1=\titleis \scriptscriptfont1=\titleiss
\textfont2=\titlesy \scriptfont2=\titlesys \scriptscriptfont2=\titlesyss
\textfont\itfam=\titleit \def\it{\fam\itfam\titleit}\rm}
 \ifx\answ\bigans\else scaled\magstep1\fi
\ifx\answ\bigans\else

 \font\absi=cmmi10 scaled\magstep1
\font\absis=cmmi7 scaled\magstep1 \font\absiss=cmmi5 scaled\magstep1
\font\abssy=cmsy10 scaled\magstep1 \font\abssys=cmsy7 scaled\magstep1
\font\abssyss=cmsy5 scaled\magstep1 
\skewchar\absi='177 \skewchar\absis='177 \skewchar\absiss='177
\skewchar\abssy='60 \skewchar\abssys='60 \skewchar\abssyss='60
\fi
\font\ninerm=cmr9 \font\sixrm=cmr6 \font\ninei=cmmi9 \font\sixi=cmmi6
\font\ninesy=cmsy9 \font\sixsy=cmsy6 \font\ninebf=cmbx9
\font\nineit=cmti9 \font\ninesl=cmsl9 \skewchar\ninei='177
\skewchar\sixi='177 \skewchar\ninesy='60 \skewchar\sixsy='60
\def\ninepoint{\def\rm{\fam0\ninerm}
\textfont0=\ninerm \scriptfont0=\sixrm \scriptscriptfont0=\fiverm
\textfont1=\ninei \scriptfont1=\sixi \scriptscriptfont1=\fivei
\textfont2=\ninesy \scriptfont2=\sixsy \scriptscriptfont2=\fivesy
\textfont\itfam=\ninei \def\it{\fam\itfam\nineit}\def\sl{\fam\slfam\ninesl}%
\textfont\bffam=\ninebf \def\bf{\fam\bffam\ninebf}\rm}
%
%

\hyphenation{anom-aly anom-alies coun-ter-term coun-ter-terms}
\def\inv{^{\raise.15ex\hbox{${\scriptscriptstyle -}$}\kern-.05em 1}}

\def\Dsl{\,\raise.15ex\hbox{/}\mkern-13.5mu D} 
\def\dsl{\raise.15ex\hbox{/}\kern-.57em\partial}

\def\lspace{\ifx\answ\bigans{}\else\qquad\fi}
\def\lbspace{\ifx\answ\bigans{}\else\hskip-.2in\fi} 
\def\boxeqn#1{\vcenter{\vbox{\hrule\hbox{\vrule\kern3pt\vbox{\kern3pt
    \hbox{${\displaystyle #1}$}\kern3pt}\kern3pt\vrule}\hrule}}}
\def\mbox#1#2{\vcenter{\hrule \hbox{\vrule height#2in
        \kern#1in \vrule} \hrule}}  
%

\def\darr#1{\raise1.5ex\hbox{$\leftrightarrow$}\mkern-16.5mu #1}

\def\roughly#1{\raise.3ex\hbox{$#1$\kern-.75em\lower1ex\hbox{$\sim$}}}

%
%


\def\frac#1#2{{#1\over#2}}

\def\journal#1&#2(#3){\unskip, #1~\bf #2 \rm(19#3) }
\def\andjournal#1&#2(#3){\sl #1~\bf #2 \rm (19#3) }

\def\bra#1{\left\langle #1\right|}
\def\ket#1{\left| #1\right\rangle}

\catcode`\@=11\def\slash#1{\mathord{\mathpalette\c@ncel{#1}}}
\overfullrule=0pt
\def\steepslash{\c@ncel}
\def\frac#1#2{{#1\over #2}}

\def\:{\!:\!}
\def\inbar{\,\vrule height1.5ex width.4pt depth0pt}
\def\IQ{\relax\,\hbox{$\inbar\kern-.3em{\rm Q}$}}
\def\IB{\relax{\rm I\kern-.18em B}}
\def\IC{\relax\hbox{$\inbar\kern-.3em{\rm C}$}}
\def\IP{\relax{\rm I\kern-.18em P}}
\def\IR{\relax{\rm I\kern-.18em R}}
\def\ZZ{\relax\ifmmode\mathchoice
{\hbox{Z\kern-.4em Z}}{\hbox{Z\kern-.4em Z}}
{\lower.9pt\hbox{Z\kern-.4em Z}}
{\lower1.2pt\hbox{Z\kern-.4em Z}}\else{Z\kern-.4em Z}\fi}

\catcode`\@=12

\def\npb#1(#2)#3{{ Nucl. Phys. }{B#1} (#2) #3}
\def\plb#1(#2)#3{{ Phys. Lett. }{#1B} (#2) #3}
\def\pla#1(#2)#3{{ Phys. Lett. }{#1A} (#2) #3}
\def\prl#1(#2)#3{{ Phys. Rev. Lett. }{#1} (#2) #3}
\def\mpla#1(#2)#3{{ Mod. Phys. Lett. }{A#1} (#2) #3}
\def\ijmpa#1(#2)#3{{ Int. J. Mod. Phys. }{A#1} (#2) #3}
\def\cmp#1(#2)#3{{ Comm. Math. Phys. }{#1} (#2) #3}
\def\cqg#1(#2)#3{{ Class. Quantum Grav. }{#1} (#2) #3}
\def\jmp#1(#2)#3{{ J. Math. Phys. }{#1} (#2) #3}
\def\anp#1(#2)#3{{ Ann. Phys. }{#1} (#2) #3}
\def\prd#1(#2)#3{{ Phys. Rev. } {D{#1}} (#2) #3}
\def\ptp#1(#2)#3{{ Progr. Theor. Phys. }{#1} (#2) #3}
\def\aom#1(#2)#3{{ Ann. Math. }{#1} (#2) #3}

\def\bs{\bigskip}

\def\bra{\langle}
\def\ket{\rangle}

\def\C{{\bf C}}

\def\N{{\bf N}}

\def\Z{{\bf Z}}
\def\cA{{\cal A}}
\def\cB{{\cal B}}

\def\cD{{\cal D}}
\def\cE{{\cal E}}

\def\cH{{\cal H}}
\def\cI{{\cal I}}

\def\cL{{\cal L}}
\def\cM{{\cal M}}

\def\cO{{\cal O}}
\def\cP{{\cal P}}

\def\cR{{\cal R}}
\def\cS{{\cal S}}

\def\cV{{\cal V}}
\def\cW{{\cal W}}

\input amssym

\def\gg{{\goth g}}

\def\gl{{\goth l}}

\def\cicy#1(#2|#3)#4{\left(\matrix{#2}\right|\!\!
                     \left|\matrix{#3}\right)^{{#4}}_{#1}}

\def\ra{\rightarrow}

\def\bs{\bigskip}

\def\Box{{\,\lower0.9pt\vbox{\hrule
\hbox{\vrule height 0.2 cm \hskip 0.2 cm
\vrule height 0.2 cm}\hrule}\,}}

\global\newcount\thmno \global\thmno=0
\def\definition#1{\global\advance\thmno by1
\bigskip\noindent{\bf Definition \secsym\the\thmno. }{\it #1}
\par\nobreak\medskip\nobreak}
\def\question#1{\global\advance\thmno by1
\bigskip\noindent{\bf Question \secsym\the\thmno. }{\it #1}
\par\nobreak\medskip\nobreak}
\def\theorem#1{\global\advance\thmno by1
\bigskip\noindent{\bf Theorem \secsym\the\thmno. }{\it #1}
\par\nobreak\medskip\nobreak}
\def\proposition#1{\global\advance\thmno by1
\bigskip\noindent{\bf Proposition \secsym\the\thmno. }{\it #1}
\par\nobreak\medskip\nobreak}
\def\corollary#1{\global\advance\thmno by1
\bigskip\noindent{\bf Corollary \secsym\the\thmno. }{\it #1}
\par\nobreak\medskip\nobreak}
\def\lemma#1{\global\advance\thmno by1
\bigskip\noindent{\bf Lemma \secsym\the\thmno. }{\it #1}
\par\nobreak\medskip\nobreak}
\def\conjecture#1{\global\advance\thmno by1
\bigskip\noindent{\bf Conjecture \secsym\the\thmno. }{\it #1}
\par\nobreak\medskip\nobreak}
\def\exercise#1{\global\advance\thmno by1
\bigskip\noindent{\bf Exercise \secsym\the\thmno. }{\it #1}
\par\nobreak\medskip\nobreak}
\def\remark#1{\global\advance\thmno by1
\bigskip\noindent{\bf Remark \secsym\the\thmno. }{\it #1}
\par\nobreak\medskip\nobreak}
\def\problem#1{\global\advance\thmno by1
\bigskip\noindent{\bf Problem \secsym\the\thmno. }{\it #1}
\par\nobreak\medskip\nobreak}
\def\others#1#2{\global\advance\thmno by1
\bigskip\noindent{\bf #1 \secsym\the\thmno. }{\it #2}
\par\nobreak\medskip\nobreak}
\def\proof{\noindent Proof: }

\def\thmlab#1{\xdef #1{\secsym\the\thmno}\writedef{#1\leftbracket#1}\wrlabeL{#1=#1}}
%
%
\def\newsec#1{\global\advance\secno by1\message{(\the\secno. #1)}
\global\subsecno=0\thmno=0\eqnres@t\noindent{\bf\the\secno. #1}
\writetoca{{\secsym} {#1}}\par\nobreak\medskip\nobreak}
\def\eqnres@t{\xdef\secsym{\the\secno.}\global\meqno=1\bigbreak\bigskip}
\def\sequentialequations{\def\eqnres@t{\bigbreak}}\xdef\secsym{}
%

%
\newcount{\exnum}
\def\prob{\advance\exnum by 1
\bigskip\item{\the\exnum.}\ }
\newcount{\exnum}
\def\next{\advance\exnum by 1
\bigskip\noindent{\the\exnum.}\ }